\documentclass[11pt]{article}
\usepackage{latexsym,amsmath,amscd,amssymb,graphics}
\usepackage{times}

\def\cF{\mathcal{F}}
\def\eps{\varepsilon}

\newtheorem{theorem}{Theorem}
\newtheorem{definition1}[theorem]{Definition}
\newtheorem{lemma}[theorem]{Lemma}
\newtheorem{remark1}[theorem]{Remark}
\newtheorem{remarks1}[theorem]{Remarks}
\newtheorem{proposition}[theorem]{Proposition}

\newtheorem{example1}[theorem]{Example}

\newenvironment{example}{\begin{example1}\normalfont}{\end{example1}}

\newenvironment{remarks}{\begin{remarks1}\normalfont}{\end{remarks1}}

\newenvironment{proof}{\noindent\textit{Proof}\quad}{\hfill
  $\blacksquare$\bigskip}

\newcommand{\bfi}{\bfseries\itshape}

\newsavebox{\savepar}

\begin{document}

\makeatletter
\title{{\bf The relation between local and global dual pairs}}
\author{James Montaldi$^{1}$, Juan-Pablo Ortega$^{2}$, and
Tudor S. Ratiu$^{3}$} \addtocounter{footnote}{1}
\footnotetext{Department of Mathematics. UMIST. PO Box 88,
Manchester M60 1QD, United Kingdom.
\texttt{j.montaldi@umist.ac.uk}. } \addtocounter{footnote}{1}
\footnotetext{Centre National de la Recherche
Scientifique. D\'epartement de Math\'ematiques de Besan\c{c}on,
Universit\'e de Franche-Comt\'e, UFR des Sciences et Techniques, 16
route de Gray, 25030 Besan\c{c}on c\'edex, France. \texttt{
Juan-Pablo.Ortega@math.univ-fcomte.fr}. }
\addtocounter{footnote}{1} \footnotetext{Centre  Bernoulli,
\'Ecole Polytechnique F\'ed\'erale de Lausanne, CH-1015 Lausanne,
Switzerland. \texttt{Tudor.Ratiu@epfl.ch}.}
\date{April 19, 2004}
\makeatother

\maketitle

\begin{abstract}
In this note we clarify the relationship between the local and
global definitions of dual pairs in Poisson geometry.  It turns
out that these are not equivalent.  For the passage from local to
global one needs a connected fiber hypothesis (this is well
known), while the converse requires a dimension condition (which
appears not to be known).  We also provide examples illustrating
the necessity of the extra conditions.
\end{abstract}

\section{Regular dual pairs}
\label{sec:regular}

The set-up we consider is the following. Let $(M, \omega)$ be a
symplectic manifold (we assume our manifolds to be paracompact),
$(P_1, \{\cdot, \cdot \}_1)$ and $(P_2, \{\cdot, \cdot \}_2)$ two
Poisson manifolds, and $\pi _1:M \rightarrow P _1$ and $\pi _2:M
\rightarrow P_2$ two surjective submersive Poisson maps. In
Remark \ref{rmk:subversive} we describe the effect of replacing
the condition ``submersion" by ``open".

Let $\cF_j$ denote the algebra of the pull-backs of smooth functions on $P_j$,  that
is,
$$\cF_j = \pi_j^*(C^\infty(P_j)).$$
Since $\pi_1$ and $\pi_2$ are Poisson it follows that  $\cF_1$
and $\cF_2$ are Poisson subalgebras of $C^{\infty}(M)$.  If $U\subset
M$ is open, we write $\cF_j(U)$ for the algebra
$$\cF_j(U) = \pi_j^*(C^\infty(\pi_j(U))).$$
This is a Poisson subalgebra of $C^\infty(U)$.

For a subset $A\subset C^\infty(U)$ we write $A^c$ for the centralizer
of $A$ with respect to the Poisson structure on $(U, \omega| _U)$,
that is,
$$A^c := \{ f\in C^\infty(U) \mid \{f, g\}^U=0 \text{ for all } g\in
A\},$$
where $\{\cdot ,\cdot \}^U$ is the restriction of the Poisson bracket to
$U$.

Note that in the following two definitions $\pi_1$ and $\pi_2$ are assumed to be
Poisson maps but not necessarily submersions.

\begin{definition1}
\label{Regular dual pairs definitions}
Consider the diagram
     \unitlength=5mm
     \begin{center}
     \begin{picture}(9,5.5)
     \put(4.9,5){\makebox(0,0){$(M, \omega)$}}
     \put(1,0.5){\makebox(0,0){$(P_1, \{\cdot, \cdot \}_1)$}}
     \put(9,0.5){\makebox(0,0){$(P_2, \{\cdot, \cdot \}_2)$}}
     \put(5,4.5){\vector(1,-1){3}}
     \put(4.5,4.5){\vector(-1,-1){3}}
     \put(2,3.5){\makebox(0,0){$\pi _1$}}
     \put(7.3,3.5){\makebox(0,0){$\pi_2$}}
     \end{picture}
     \end{center}
\begin{itemize}
\item The diagram forms a {\bfi Howe (H) dual pair\/} if the Poisson
subalgebras $\cF_1$ and $\cF_2$ centralize each other:
\begin{equation}
\label{double commutant relation expression}
\cF_1^c = \cF_2\quad\mbox{and}\quad \cF_2^c = \cF_1.
\end{equation}

\item The diagram forms a {\bfi Lie-Weinstein (LW) dual pair\/}
when $\ker T \pi_1$ and $\ker T \pi_2$ are symplectically orthogonal
distributions. That is, for each $m \in M $,
\begin{equation}
\label{condition lie weinstein point}
\left( \ker T _m \pi_1 \right) ^\omega= \ker T _m \pi_2.
\end{equation}
\end{itemize}
In each case, the dual pair is {\bfi regular} when the maps
$\pi_j$ are assumed to be surjective submersions; otherwise they
are {\bfi singular} dual pairs.
\end{definition1}

This notion of singular dual pair is less general than that
in~\cite{ortega02}. Notice that for a regular Lie-Weinstein
dual pair, the dimensions of $P_1$ and $P_2$ sum to the dimension
of $M$. Actually, we note that if the manifold $M$ is Lindel\"of
or paracompact as a topological space then the Lie-Weinstein
condition cannot hold unless the dual pair is regular.  This is
because the LW condition implies that the two maps $\pi  _1 $ and
$\pi _2 $ are of complementary rank (the ranks sum to $\dim(M)$)
and, by the lower semicontinuity of the rank of a smooth map, the
maps must both be of constant rank. Since they are surjective they
must be submersions (by Sard's theorem).

We emphasize that the definition of a Lie-Weinstein dual pair is
local, while that of a Howe dual pair is global.  However, the
latter definition can be localized as follows.

\begin{definition1}
The diagram above forms a {\bfi local Howe (LH) dual pair\/} if for each $m\in M$
and each neighbourhood $V$ of $m$ there is a neighbourhood $U$ of $m$ with $U\subset
V$ such that the algebras $\cF_1(U)$ and $\cF_2(U)$ centralize each other in
$C^\infty(U)$. If in addition, $\pi_1$ and $\pi_2$ are surjective submersions, then
the local Howe dual pair is said to be {\bfi regular\/}; otherwise it is {\bfi
singular\/}.
\end{definition1}

The notion of Howe dual pair has its origins in the study of group representations
arising in quantum mechanics (see for instance~\cite{howe, kashiwara vergne,
Sternberg Wolf 1978, jakobsen vergne}, and references therein) and it appears for
the first time in the context of Poisson geometry in~\cite{weinstein83}.  The
definition of Lie-Weinstein dual pair can be traced back to~\cite{lie1890} and, in
its modern formulation, is due to~\cite{weinstein83}.  Examples of dual pairs
arising in classical mechanics can be found in~\cite{semidirect marsden ratiu
weinstein, mwr clebsch}, and references therein. Further details on dual pairs can
also be found in \cite{the book!}.

The relationships between the three notions of regular dual pair can
be summed up in the following two results.

\begin{proposition}
\label{prop:local}
The two local notions of regular dual pair, that is, Lie-Weinstein
and local Howe, are equivalent.
\end{proposition}

In Remark \ref{rmk:subversive} we provide an example showing that
this result no longer holds if the regularity is dropped.

\begin{theorem}\
\label{thm:local-global}
\begin{enumerate}
\item
If a regular Howe dual pair is such that the Poisson manifolds $P_1$ and $P_2$ are
of complementary dimension, that is, $\dim P_1+\dim P_2=\dim M$, then it forms a
regular local Howe dual pair.
\item
If a regular local Howe dual pair is such that the fibers of $\pi_1$ and $\pi_2$ are
connected then it is a regular Howe dual pair.
\end{enumerate}
\end{theorem}

Before giving the proofs (in Section \ref{sec:proofs}), we give
two examples showing the necessity of the hypotheses.

\begin{example}
\label{ex:LH not H}
This example (suggested to us by Andrea Giacobbe) shows that the
hypothesis of connected fibers in the passage from local to global
in the theorem above is necessary.  Let
$$
\mathbb{T}^2 =\left\{ (\theta_1,\theta_2)\mid \theta _j\in
\mathbb{R}/2\pi\mathbb{Z} \right\}
$$
be the 2-torus considered as a symplectic manifold with the area form $\omega:=
\mathbf{d} \theta _1\wedge \mathbf{d} \theta _2 $. Consider the diagram
$S^1\stackrel{\pi_1}{\leftarrow} \mathbb{T}^2\stackrel{\pi_2}{\rightarrow} S^1$ with
$ \pi _j(\theta_1, \theta _2):= j\theta_1 $. The fibers of $\pi _2 $ have two
connected components. It is easy to see that this forms a Lie-Weinstein dual pair
(and hence a local Howe dual pair) but not a Howe dual pair.  Indeed, the function
$\cos (\theta _1) $ belongs to $\cF_1^c $ but not to $\cF_2$.
\end{example}

In the example below and in subsequent proofs, we use the
following notation. On the symplectic manifold $(M,\omega)$, we
write the Poisson tensor as $B\in\Lambda^2(M)$. If $h\in
C^\infty(M)$ then the Hamiltonian vector field $X_h$ is defined by
$\mathbf{d}h =\mathbf{i}_{X_h}\omega:= \omega(X_h,\cdot)$. The Poisson
tensor, defined by $B(\mathbf{d}g, \mathbf{d}h) := X_h[g] = \langle
\mathbf{d}g, X_h \rangle = \omega(X_g,\,X_h)$, induces the vector bundle
morphism over the identity $B^\sharp:T^*M\to TM$ given by
$B^\sharp(\mathbf{d}h) = X_h$.

\begin{example} \label{ex:Howe not L-W}
This example shows that without the dimension hypothesis, a regular Howe dual pair
need not be locally Howe (nor Lie-Weinstein).  Let $M:= \mathbb{T} ^3 \times
\mathbb{R} $ and $\lambda _1, \lambda _2, \lambda _3 \in \mathbb{R} $ be linearly
independent over $\mathbb{Q}$. Define on $M$ the symplectic structure $\omega$ whose
Poisson tensor $B \in \Lambda^2(M)$ is given by
\[
B ^{\sharp}= \left(
\begin{array}{cccc}
           0 &1 &0 &\lambda _1\\
           -1 &0 &0  & \lambda _2\\
           0 &0 &0 &\lambda _3\\
           - \lambda _1&- \lambda _2&- \lambda _3 &0
\end{array}
\right).
\]
Let $\pi:\mathbb{T} ^3 \times \mathbb{R} \rightarrow \mathbb{R} $ be
the projection onto the $\mathbb{R}$ factor. The Hamiltonian vector
field associated to the function $\pi$ is
$$X_\pi = \lambda_1\frac{\partial}{\partial\theta_1} +
\lambda_2\frac{\partial}{\partial\theta_2} +
\lambda_3\frac{\partial}{\partial\theta_3}.$$ Then the diagram
$\mathbb{R}\stackrel{\pi }{\leftarrow} M \stackrel{\pi}{\rightarrow} \mathbb{R}$ is
a regular Howe pair but clearly not a Lie-Weinstein dual pair and hence not a
regular local Howe dual pair.  In order to see that it is a regular Howe dual pair
let $g \in \left( \pi ^\ast C^{\infty} (\mathbb{R})\right) ^c $.  The trajectories
of the vector field $X _\pi $ on $M$ are irrational windings which are dense in the
fibers of $\pi $. Since $g$ is invariant under this Hamiltonian flow it must be
constant on the fibers of $\pi$ and hence $g \in \pi^\ast C^{\infty}(\mathbb{R}) $,
as required.
\end{example}

\begin{remarks}
\label{rmk:subversive}\

\textbf{1.}
If one merely assumes the maps $\pi_j$ to be open rather than
submersions, then one can still pass from local Howe to global, provided of
course the fibers are connected (see the proof in Section
\ref{sec:proofs}).

\textbf{2.}  As was already observed, if the Lie-Weinstein
condition is satisfied then the two Poisson maps $\pi_j$ are of
constant rank. Any example of a local Howe dual pair which is not
Lie-Weinstein must of course be singular (the maps cannot be
submersions) and not of constant rank. A simple example is
obtained by putting $M=\mathbb{R}^4$ with coordinates
$(x_1,x_2,y_1,y_2)$ and its usual symplectic form $\omega=\sum_j
\mathbf{d} x_j\wedge \mathbf{d} y_j$. Let $\pi_1:M\to\mathbb{R}$
be given by $\pi_1(\mathbf{x},\mathbf{y}) =
\mathbf{x}\cdot\mathbf{y}$ (inner product) and $\pi_2:M\to
\mathbb{R}^4$ be given by $\pi_2(\mathbf{x}, \mathbf{y}) =
\mathbf{x}\otimes\mathbf{y}$ (outer product!).  In coordinates,
$$\pi_2(x_1,x_2,y_1,y_2) = (x_1y_1,x_1y_2,x_2y_1,x_2y_2).$$
The fibers of both $\pi_1$ and $\pi_2$ are connected and the image
$\pi_2(M)$ is a cone in $\mathbb{R}^4$.

We claim that $\mathbb{R}\stackrel{\pi_1}{\leftarrow} M
\stackrel{\pi_2}{\rightarrow}\mathbb{R}^4$ forms a singular Howe
dual pair.  It is clearly not Lie-Weinstein at the origin, as
$T\pi_1$ and $T\pi_2$ both vanish at that point.  That it is a Howe
pair follows from the paper of Karshon and Lerman \cite{KL97},
since $\pi_1$ is the orbit map for the $U(2)$ action on
$\mathbb{C}^2$ and $\pi_2$ its momentum map.
Note that in \cite{montaldi tokieda} it is shown that momentum maps of
representations are $G$-open, although we do not know whether $G$-openness is
sufficient to be able to pass from local to global in the singular case.
\end{remarks}

\begin{figure}[htb]
\begin{center}
\includegraphics{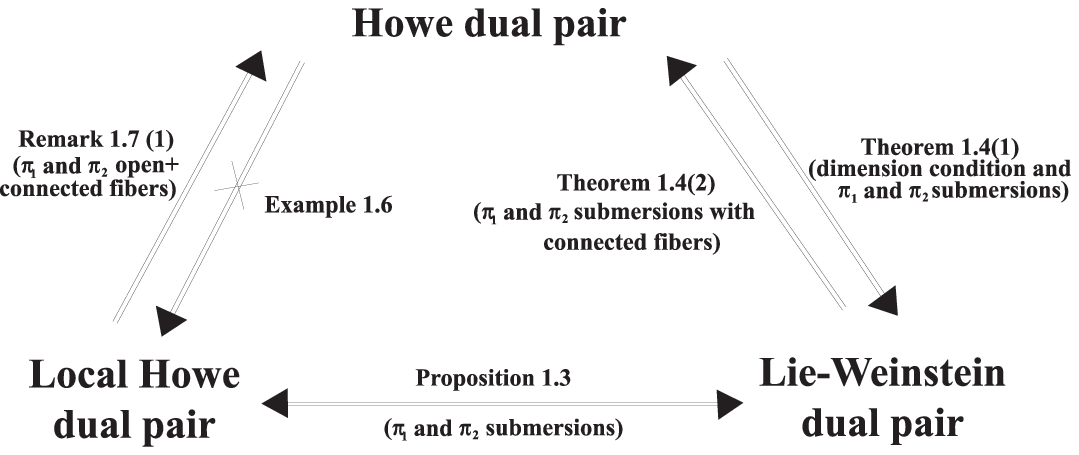}
\end{center}
\end{figure}


\section{Proofs}
\label{sec:proofs}

Define $K_j:=\ker T\pi_j$ (for $j=1,2$), which are two
subbundles of $TM$. Since the maps $\pi_j$ are submersions, we
have for each $m\in M$
\begin{equation}
\label{eq:annihilator}
K_j(m)^\circ=\{\mathbf{d}f(m) \mid f\in\cF_j\}.
\end{equation}
Furthermore, since $K_j(m)^\omega=B^\sharp(K_j(m)^\circ)$ it follows that
\begin{equation}
\label{eq:Komega}
K_j(m)^\omega = \{X_f(m) \mid f\in\cF_j\}.
\end{equation}

\subsection*{Proof of Proposition \ref{prop:local}}
We establish the equivalence of the two local notions of regular dual
pair.  Recall that  $\pi_1$ and $\pi_2$ are assumed to be surjective
submersions.

\paragraph{LH $\Rightarrow$ LW:}
We wish to show that $(K_1)^\omega=K_2$, which we do by double inclusion.

First we show $(K_1)^\omega\subset K_2$.  Let $z_j$ be coordinates on
$P_2$. Then by the local Howe condition, $f\in\cF_1$ implies
$$
0 = \{f,\,z_j\circ\pi_2\} = \omega(X_f, X_{z_j \circ \pi_2} )\quad
\text{for all} \quad j
$$
so that $X_f\in \{ X_g \mid g \in \mathcal{F}_2 \} ^ \omega = K_2$,
by (\ref{eq:Komega}).

For the converse inclusion, $K_2\subset (K_1)^\omega$, let $m\in
M$ and $v\in K_2(m)$.  Then there is a function $f$ such that
$v=X_f(m)$. Then for $g\in \cF_2$ we have
$\left<\mathbf{d}f(m),\,X_g(m)\right> =-\mathbf{d}g(m)(v)=0$, so
that $\mathbf{d}f(m) \in (K_2(m)^\omega)^\circ$. Since
$(K_2{}^\omega)^\circ$ is a subbundle of $T^*M$ and hence locally
trivial, we can choose $f$ on a neighbourhood $V$ of $m$ so that
at each $x\in V$, $\mathbf{d}f(x) \in (K_2(x)^\omega)^\circ$.
Thus
$$
0 = \left<\mathbf{d}(z_j\circ\pi_2)(x),X_f(x) \right> =
\{z_j\circ\pi_2,f\}(x).
$$
By the local Howe condition, this implies there is a sub-neighbourhood $U\subset V$
of $m$ such that $f\in \cF_1(U)$. Consequently, $X_f$ is a section of $
(K_1)^\omega$ over $U $ and, since $v=X_f(m)$, it follows that
$K_2\subset(K_1)^\omega$.

\paragraph{LW $\Rightarrow$ LH:}
We choose $U$ such that both submersions $\pi_j|_U$ have connected fibers (this can
be done using Lemma \ref{lem:connected}, by choosing any function $h$ on $M$ with a
non-degenerate local minimum at $m$). We prove the equality $\cF_2(U) = \cF_1(U)
^c$, again by double inclusion.

Let $f \in \cF_2(U)$.  Since  $ \mathbf{d}f(m) \in K_2(m)^{\circ}$, for any
$m\in  U$ and, by hypothesis, $(K_1) ^\omega= K_2$, we
have
\begin{equation*}
X_{f}(m) \in B ^{\sharp} (m) (K_2(m)^{\circ})  = K_2(m)^\omega =
K_1(m)=\left(\{\mathbf{d}g(m)\mid g \in \cF_1(U)\}\right)^{\circ},
\end{equation*}
where the last equality follows from~(\ref{eq:annihilator}).
Consequently, for an arbitrary $g \in \cF_1(U)$, we conclude
\[
\{ g , f\}= \mathbf{d}g\left(X_{f} \right) =0,
\]
which implies that $f \in \mathcal{F}_1(U)^c $.

Conversely, let $f \in \mathcal{F} _1 (U)^c $ for some $U$. In
order to prove that $f \in \mathcal{F}_2(U)$ we start by showing
that it is locally constant on the fibers of $\pi _2 $. Indeed,
since $\pi_2 $ is a surjective submersion and the diagram $(P_1,
\{\cdot, \cdot\}_1)\stackrel{\pi_1}{\leftarrow} (M,
\omega)\stackrel{\pi_2}{\rightarrow} (P _2,\{\cdot, \cdot\}_2)$
forms a Lie-Weinstein dual pair, for any $m \in M $ we have
\[
T _m \left(\pi_2 ^{-1}( \pi_2 (m)) \right)= K_2(m)=
\left(K_1(m)\right)^\omega= B ^{\sharp} (m)
\left((K_1(m))^{\circ}\right).
\]
This equality, together with~(\ref{eq:annihilator}), guarantees that
any vector $v$ tangent at $m$ to the fiber $\pi_2 ^{-1}( \pi_2 (m)) $
can be written as $v=X_g(m) $, for some $g \in \cF_1(U)$. Hence,
\[
\mathbf{d} f(m) (v) =\mathbf{d} f(m)\left(X_{g}(m)\right) =
\{f, g\} (m)=0.
\]
Since both $m\in M$ and $v \in T _m \left(\pi_2 ^{-1}( \pi_2 (m)) \right)$ are
arbitrary and, by the choice of $U$, the fibers of $\pi_2 | _U$ are connected, this
equality guarantees that the function $f \in \cF_1(U)^c$ is constant on the fibers
of $\pi_2 $.  This implies that there exists a unique function $ \overline{f} :\pi
_2 (U)\rightarrow \mathbb{R}$ that satisfies the equality $f= \overline{f} \circ
\pi_2|_U$. Since $f$ is smooth and $\pi_2|_U$ a submersion it follows that
$\overline{f} $ is smooth and hence $f \in \cF_2(U)$.

The equality $\cF_1(U) =\cF_2(U)^c $ is proved analogously.\hfill
     $\blacksquare$

\subsection*{Proof of Theorem \ref{thm:local-global}}

\paragraph{Global $\Rightarrow$ Local:}
Let $m \in M $ be an arbitrary point. We will now show that the
hypotheses in the statement imply that $(K_1(m))^\omega =K_2(m)$.
The  inclusion $(K_2 (m)) ^\omega \subset K_1 (m) $ follows from the
global Howe hypothesis and the equalities~(\ref{eq:Komega})
and~(\ref{eq:annihilator}). Indeed,
\[
(K_2 (m)) ^\omega=\{X _f (m)\mid f \in \mathcal{F} _2\}
=\{X _g(m)\mid g \in \mathcal{F} _1 ^c\}\subset
\left(\{ \mathbf{d} f(m)\mid f \in \mathcal{F} _1\}
\right)^{\circ}=K _1 (m).
\]
The converse inclusion follows immediately from the hypothesis on
the dimensions and the submersiveness of $\pi_1$ and $\pi_2 $.
Indeed,
\[
\dim (K_2(m)) ^\omega=\dim P _2=\dim M-\dim P _1=\dim K_1(m),
\]
and hence $(K_2(m)) ^\omega = K_1(m) $, as required.


\paragraph{Local $\Rightarrow$ Global:}
Using the symmetry of the statement with respect to the exchange
of $\pi _1 $ with $\pi _2 $ it suffices to show that, for
instance, $\cF_1^c = \cF_2$.  Notice that the proof below only
requires the Poisson maps $\pi_j$ to be open onto $P_j$, rather
than submersions, and that the $P_j$ are not even required to be
manifolds (in which case the proof would hold with an appropriate
algebraic definition of smooth function on $P _j$).

First, it is easy to show that given $f \in\cF_1$ and $g\in\cF_2$
then $\{f, g\} =0 $.  Indeed, let $m\in M$ and let $U$ be a
neighbourhood of $m$ on which the local Howe condition holds.
Then, $f| _{U} \in \cF_1(U)$ and similarly $g|_{U} \in \cF_2(U)$
and hence,
\begin{equation}
\label{needed relation}
\{f, g\} (m)=\{f|_{U}, g|_{U}\}^{U} (m)=0,
\end{equation}
where $\{\cdot , \cdot \}^U$ is the restriction of the Poisson bracket to
$U$. Since $m\in M$, $f \in \cF_1$, and $g \in\cF_2$ are arbitrary, it
follows that $\cF_2\subset\cF_1^c.$

Second, let $g \in \cF_1^c $; one needs to show that $g \in
\cF_2$. For each $m\in M$ there is a neighbourhood $U_m$ on which
the local Howe hypothesis is valid.  This provides a cover of $M$,
from which we can extract a locally finite subcover $\{U_a\}_{a\in
A}$.

For each $U_a$, we claim that $g| _{U_a} \in \cF_1(U_a)^c$.  It then
follows by hypothesis that $g| _{U_a} \in \cF_2(U_a)$, which allows us
to write
\begin{equation}
\label{pull back gives us another}
g|_{U_a}= \overline{g}_a\circ  \pi _2|_{U_a}
\end{equation}
for some $\overline{g} _a \in C^{\infty}(\pi _2 (U_a))$. Our second
claim is that there is a function $\overline{g} \in C^{\infty}(P _2) $
such that $\overline{g} _a= \overline{g}|_{\pi _2 (U _a)} $. The result
then follows as $g= \overline{g} \circ \pi _2$. We now establish the two
claims.

The first claim is proved by contradiction. Suppose that there
exists a function $f \in C^{\infty}(\pi _1 (U _a)) $ and $x \in U
_a $ such that $\{g|_{U _a}, f \circ \pi_1|_{U _a}\}^{U _a} (x)
\neq 0 $.  Let $V _x $ be an open neighbourhood of $x$ such that
$\overline{V _x} \subset \pi _1(U _a)$. Then there is an extension
$F \in C^{\infty}(P_1)$ of $f|_{V _x}$. Since $g \in \cF_1^c$ it
follows that
\[
0=\{g, F \circ \pi _1\} (x)=\{g|_{U _a}, f \circ \pi _1|_{U _a}\} ^{U_a}
(x)\neq 0,
\]
contradicting \eqref{needed relation}.

As to the second claim, notice first that~(\ref{pull back gives us
another}) implies that $g$ is locally constant along the fibers of
$\pi _2 $. Since by hypothesis these fibers are connected, $g$ is
constant on the fibers of $\pi _2 $ and $\overline{g} $ is therefore
well defined. Moreover, $\overline{g} $ coincides with $\overline{g}_a $
on the open sets of the form $\pi_2(U _a)$ and so is smooth.
\hfill $\blacksquare$

\begin{lemma}
\label{lem:connected} Let $U$ be a manifold and $h$ be a smooth real-valued function
with a non-degenerate local minimum at $u_0\in U$, and suppose $h(U_0)=0$. Let
$\pi:U\to P$ be a submersion in a neighbourhood of $u_0$. Then for $\eps$
sufficiently small the level sets of $\pi$ restricted to $B_\eps$ are diffeomorphic
to an open ball, where $B_\eps$ is the connected component of $\{u\in U\mid
h(u)<\eps\}$ containing $u_0$.
%
%
\end{lemma}

\begin{proof}
Write $F_{y,\eps} = B_\eps\cap\pi^{-1}(y)$. We wish to show that when it is
nonempty, $F_{y,\eps}$ is diffeomorphic to an open ball.  First, choose coordinates
near $u_0$ and $\pi(u_0)$ such that $\pi(x,y)=y$. Since the restriction of $h$ to
$\{y=0\}$ has a non-degenerate minimum at $x=0$, by the splitting lemma (or Morse
lemma with parameters, see e.g.\ \cite[p.~97]{BG}) there is a neighbourhood $U_1$ of
$u_0$ on which one can change coordinates by $(x,y)\mapsto(X,y)=(X(x,y),y)$ in such
a way that $h(X,y) = Q(X)+g(y)$, where $Q$ is a positive definite quadratic form and
$g$ is a smooth function. Choose $\eps_1$ sufficiently small so that $B_{\eps_1}$ is
contained in this neighbourhood. In these coordinates, for $\eps\leq\eps_1$ and for
each $y\in \pi(U_1)$,
$$F_{y,\eps}=\{(X,y)\in B_\eps\mid Q(X)<\eps-g(y)\}.$$
It is clear that for each $y,\eps$ this set is either empty or diffeomorphic to an
open ball.
\end{proof}


\section{Final remarks}
\label{sec:rmks}

There are a number of theorems in the literature which are stated with the
hypothesis that a given set-up is a Howe dual pair, while the proof uses the
Lie-Weinstein property. The most famous of these are probably the Symplectic Leaves
Correspondence Theorem~\cite{weinstein83, blaom dual pairs} and Weinstein's theorem
on transverse Poisson structures \cite[Theorem 8.1]{weinstein83}.  Here we give an
example showing that these theorems fail if one does not assume a local
(Lie-Weinstein) hypothesis.  Of course, by the theorem above the hypothesis that
$P_1$ and $P_2$ are of complementary dimension is also sufficient.

\begin{example}
Let
\begin{equation}
\label{eq:T5xR}
M=\mathbb{T}^2 \times \mathbb{T}^3 \times \mathbb{R}
\end{equation}
with coordinates $(\theta,\phi,x)$, where
$\theta=(\theta_1,\theta_2)\in\mathbb{T}^2$,
$\phi=(\phi_1,\phi_2,\phi_3)\in\mathbb{T}^3$, and
$\theta_j,\phi_j\in\mathbb{R}/2\pi\mathbb{Z}$.  Let $\pi_1:M\to
\mathbb{T}^3 \times \mathbb{R}$ and $\pi_2:M\to \mathbb{R}$ be given
by
$$\pi_1(\theta,\phi,x) = (\phi,x),\quad \mbox{and}\quad
\pi_2(\theta,\phi,x) = x.$$
We endow $M$ with a Poisson structure whose Poisson tensor can be
written in block form as
$$B^\sharp = \left[\begin{matrix}0&A&B\cr -A^T&C&0\cr
-B^T&0&0\end{matrix}\right].$$ We assume that the entries of $A,B,C$ are independent
over $\mathbb{Q}$, excluding those forced to vanish. An argument along similar lines
to that for Example~\ref{ex:Howe not L-W} shows that
$P_1\stackrel{\pi_1}{\leftarrow} M \stackrel{\pi_2}{\rightarrow} P_2$ forms a
regular Howe dual pair. However, since the dimensions of $P_1=\mathbb{T}^3 \times
\mathbb{R}$ and $P_2=\mathbb{R}$ differ by an odd number, it is clear that they
cannot have the property of having anti-isomorphic transverse Poisson structures, up
to a product with a symplectic factor.  Indeed, the Poisson structure on $P_2$ is of
course trivial, while that on $P_1$ is of rank 2.

\medskip

This example also shows that in the absence of the local
Lie-Weinstein hypothesis, the correspondence between the
symplectic leaves of two Poisson manifolds in duality may fail.
In order to see this take $x \in \mathbb{R}$. This is a symplectic
leaf in $\mathbb{R}$, and if the Symplectic Leaves Correspondence
theorem were valid then $\pi_1(\pi_2^{-1}(x))\simeq \mathbb{T}^3$
would be a symplectic leaf in $\mathbb{T}^2\times
\mathbb{T}^3\times\mathbb{R}$.  This is obviously impossible as
$\mathbb{T}^3$ is of odd dimension.
\end{example}

\bigskip

\noindent\textbf{Acknowledgments.}
We thank Ana Cannas da Silva, Eugene Lerman, Rui Loja Fernandes, and Alan
Weinstein for patiently answering our questions and for several
illuminating discussions. This research was partially supported by the
European Commission and the Swiss Federal Government through funding for
the Research Training Network \emph{Mechanics and Symmetry in Europe}
(MASIE). The support of the Swiss National Science Foundation is also
acknowledged.

\end{document}